\documentclass[notitlepage,10pt]{article}
\catcode`\@=11
\@addtoreset{equation}{section}

\catcode`\@=12

\usepackage{latexsym}
\usepackage{amsmath}
\usepackage{amsfonts}
\usepackage{amssymb}
\usepackage{mathrsfs}
\usepackage{comment}

 \usepackage[colorlinks=true]{hyperref}
\hypersetup{urlcolor=blue, citecolor=red}

\newtheorem{Theorem}{Theorem}[section]

\newtheorem{Lemma}[Theorem]{Lemma}
\newtheorem{Proposition}[Theorem]{Proposition}

\def\R{\mathbb{R}}
\def\N{\mathbb{N}}
\def\S{\mathbb{S}}

\def\proof{\noindent{\textbf{Proof. }}}
\def\QED{\hfill {$\square$}\goodbreak \medskip}

\begin{document}

\title{Entire solutions for a class of fourth order semilinear elliptic equations with weights}

\author{}



\date{}

{\Large
\centerline{Entire solutions for a class of fourth order}
\centerline{semilinear elliptic equations with weights}
}
\bigskip

\centerline{\scshape Paolo Caldiroli\footnote{email: paolo.caldiroli@unito.it} and Gabriele Cora\footnote{email: 330933@studenti.unito.it}}
\medskip

{\footnotesize
 \centerline{Dipartimento di Matematica, Universit\`a di Torino}
 \centerline{via Carlo Alberto, 10 -- 10123 Torino, Italy}
}
\bigskip

\begin{abstract}
We investigate the problem of entire solutions for a class of fourth order, dilation invariant, semilinear elliptic equations with power-type weights and with subcritical or critical growth in the nonlinear term. These equations define non compact variational problems and are characterized by the presence of a term containing lower order derivatives, whose strength is ruled by a parameter $\lambda$. We can prove existence of entire solutions found as extremal functions for some Rellich-Sobolev type inequalities. Moreover, when the nonlinearity is suitably close to the critical one and the parameter $\lambda$ is large, symmetry breaking phenomena occur and in some cases the asymptotic behavior of radial and non radial ground states can be somehow described. 
\vspace{6pt}\\
\textbf{Keywords:} Weighted biharmonic operator, extremal functions, Rellich-So\-bo\-lev inequality, breaking symmetry.\vspace{6pt}\\
\textbf{2010 Mathematics Subject Classification:} 26D10, 47F05.
\end{abstract}

\section{Introduction}

In recent years much interest has been addressed to a class of equations shaped on 
\begin{equation}
\label{eq:NLRellich}
\Delta(|x|^{\alpha}\Delta u)=|x|^{-\beta}|u|^{q-2}u\quad\text{in }\R^{n}
\end{equation}
where the dimension $n$ and the parameters $\alpha$, $\beta$ and $q>2$ are asked to satisfy suitable restrictions. In particular, in the case of the pure biharmonic operator, we quote, e.g., \cite{AOM}, \cite{BerFarFerGaz}, \cite{B}, \cite{BM}, \cite{C14}, \cite{Lions}, \cite{Mit93}, \cite{Mit96}, \cite{Mus13}, \cite{NousSwanYang}.

Equations like (\ref{eq:NLRellich}) arise in a natural way from variational inequalities of the form
$$
S_{\alpha,q}\left(\int_{\R^{n}}|x|^{-\beta}|u|^{q}~\!dx\right)^{2/q}\le\int_{\R^{n}}|x|^{\alpha}|\Delta u|^{2}~\!dx\quad\forall u\in C^{\infty}_{c}(\R^{n}\setminus\{0\})
$$
which, for $q>2$, can be considered as nonlinear versions of the weighted Rellich inequality (this case occurs taking $q=2$ and $-\beta=\alpha-4$). 

In addition, (\ref{eq:NLRellich}) can be viewed as a higher order version of equations like
\begin{equation}
\label{eq:NLHardy}
-\mathrm{div}(|x|^{a}\nabla u)=|x|^{-b}|u|^{q-2}u\quad\text{in }\R^{n}
\end{equation}
where again $n\in\N$, $a,b\in\R$ and $q>2$ are subject to some constraints. Even (\ref{eq:NLHardy}) comes from a class of variational inequalities which are often known in the literature as Caffarelli-Kohn-Nirenberg inequalities (\cite{CKN}) and can be read as interpolation of the linear weighted Hardy inequality with the weighted Sobolev inequality.

As one can expect, both for (\ref{eq:NLRellich}) and for (\ref{eq:NLHardy}), the restrictions on the parameters are needed in order to guarantee that the supporting variational inequalities hold true. In particular, these restrictions necessarily involve some dilation invariance which is a typical feature of any problem displaying scaling processes.  

In this work we study a class of equations which are built as linear combination of (\ref{eq:NLRellich}) and (\ref{eq:NLHardy}). In particular we are interested in the existence of nontrivial solutions to 
\begin{equation}
\label{eq:pb}
\left\{\begin{array}{l}
\Delta(|x|^{\alpha}\Delta u)-\lambda~\!\mathrm{div}(|x|^{\alpha-2}\nabla u)=|x|^{-\beta}|u|^{q-2}u\quad\text{in }\R^{n}\\
\int_{\R^{n}}|x|^{\alpha}|\Delta u|^{2}~\!dx<\infty
\end{array}\right.
\end{equation}
where $n\ge 5$, $q>2$, 
\begin{equation}
\label{eq:alpha-beta}
4-n<\alpha<n~\!,~~
\beta=n-\frac{q(n-4+\alpha)}{2}~\!,
\end{equation}
and $\lambda$ is a real parameter subject to some limitation. 
In particular we look for ground states of (\ref{eq:pb}), i.e., solutions to (\ref{eq:pb}) characterized as minimizers for the following problems:
\begin{equation}
\label{eq:Sq}
S_{\alpha,q}(\lambda):=\inf_{\scriptstyle u\in D^{2,2}(\R^{n};|x|^{\alpha})\atop\scriptstyle u\ne 0}\frac{\int_{\R^{n}}|x|^{\alpha}|\Delta u|^{2}~\!dx+\lambda\int_{\R^{n}}|x|^{\alpha-2}|\nabla u|^{2}~\!dx}{\left(\int_{\R^{n}}|x|^{-\beta}|u|^{q}~\!dx\right)^{{2}/{q}}}~\!.
\end{equation}
Here $D^{2,2}(\R^{n};|x|^{\alpha})$ is the space defined as the completion of $C^{\infty}_{c}(\R^{n})$ with respect to the norm
\begin{equation}
\label{eq:norm}
\|u\|_{2,\alpha}^{2}:=\int_{\R^{n}}|x|^{\alpha}|\Delta u|^{2}~\!dx~\!.
\end{equation}
We point out that the role of entire solutions and especially of ground states of (\ref{eq:pb}) is rather meaningful since this kind of solutions naturally appear as limiting profiles in the blowup analysis of related classes of nonlinear problems.

As discussed later, thanks to already known variational inequalities (see \cite{CM-Milan11}, \cite{GM13}), one has that $S_{\alpha,q}(\lambda)>0$ and the minimization problem (\ref{eq:Sq}) makes sense whenever 
\begin{equation}
\label{eq:q}
2<q\le 2^{**}:=\frac{2n}{n-4}
\end{equation}
and 
\begin{equation}
\label{eq:delta-nabla}
\lambda>-\gamma_{\alpha}\text{~~where~~}\gamma_{\alpha}:=\inf_{\scriptstyle u\in D^{2,2}(\R^{n};|x|^{\alpha})\atop\scriptstyle u\ne 0}\frac{\int_{\R^{n}}|x|^{\alpha}|\Delta u|^{2}~\!dx}{\int_{\R^{n}}|x|^{\alpha-2}|\nabla u|^{2}~\!dx}~\!.
\end{equation}
In particular it is known that if $n\ge 5$ and $\alpha\in(4-n,n)$ then $\gamma_{\alpha}>0$. Moreover for $\alpha\in[0,n)$ then $\gamma_{\alpha}=(n-\alpha)^{2}/4$. These facts are discussed in \cite{CM-Bangalore}, \cite{Mor}, \cite{GM11}. 

Taking $\beta$ as in (\ref{eq:alpha-beta}) makes problem (\ref{eq:pb})  invariant with respect to the action of the weighted dilation group
\begin{equation}
\label{eq:dilation}
\rho\mapsto (\rho*u)(x)=\rho^{\frac{n-4+\alpha}{2}}u(\rho x)\quad(\rho>0).
\end{equation}
This invariance is responsible of a lack of compactness in the study of the minimization problem (\ref{eq:Sq}). Adapting some techniques already used for different problems (\cite{BM}, \cite{CM-Milan11}), we develop a suitable argument allowing us to recover some local compactness, and we get a first existence result, stated as follows. 

\begin{Theorem}
\label{T:ground-state}
Let $n\ge 5$ and assume (\ref{eq:alpha-beta}). Then:
\begin{itemize}
\item[(i)]
For $q\in(2,2^{**})$ and $\lambda\in(-\gamma_{\alpha},\infty)$ problem (\ref{eq:pb}) admits a ground state.
\item[(ii)]
For $q=2^{**}$, for every $\alpha\in(4-n,n)$ there exists $\lambda_{\alpha}>-\gamma_{\alpha}$ such that problem (\ref{eq:pb}) admits a ground state if $\lambda\in(-\gamma_{\alpha},\lambda_{\alpha})$ (the value $\lambda_{\alpha}$ is given by (\ref{eq:lambda-n-alpha})). 
\end{itemize}
\end{Theorem}


By exploiting the rotational symmetry of the domain and of the weights in (\ref{eq:pb}), we can drop the upper bound on $q$ in (\ref{eq:q}) by looking for radial ground states for problems (\ref{eq:pb}), namely, non trivial, radial weak solutions of (\ref{eq:pb}) characterized as minimum points for 
\begin{equation}
\label{eq:Sqrad}
S_{\alpha,q}^{\mathrm{rad}}(\lambda):=\inf_{\scriptstyle u\in D^{2,2}_{\text{rad}}(\R^{n};|x|^{\alpha})\atop\scriptstyle u\ne 0}\frac{\int_{\R^{n}}|x|^{\alpha}|\Delta u|^{2}~\!dx+\lambda\int_{\R^{n}}|x|^{\alpha-2}|\nabla u|^{2}~\!dx}{\left(\int_{\R^{n}}|x|^{-\beta}|u|^{q}~\!dx\right)^{{2}/{q}}}
\end{equation}
where $D^{2,2}_{\text{rad}}(\R^{n};|x|^{\alpha})$ is the space of radial functions belonging to $D^{2,2}(\R^{n};|x|^{\alpha})$. We have that:

\begin{Theorem}
\label{T:radial-ground-state}
Let $n\ge 5$ and assume (\ref{eq:alpha-beta}). Then for every $q\in(2,\infty)$ and $\lambda>-(n-\alpha)^{2}/{4}$, problem (\ref{eq:pb}) admits a radial ground state. Moreover such a ground state has constant sign and is unique up to the weighted dilation (\ref{eq:dilation}).
\end{Theorem}

The second part of our work consists in the study of global ground states of (\ref{eq:pb}) given by Theorem \ref{T:ground-state}. In particular we are interested in investigating radial symmetry or not of these solutions. We find symmetry breaking in different situations.
A first result in this direction is the following.

\begin{Theorem}
\label{T:SB1}
Let $n\ge 5$ and $\lambda>0$. There exist $\overline\alpha>0$ and $\overline{q}\in(2,2^{**})$, both depending on $\lambda$, such that if $q\in(\overline{q},2^{**}]$ and $|\alpha|<\overline\alpha$ then $S_{\alpha,q}(\lambda)<S_{\alpha,q}^{\mathrm{rad}}(\lambda)$. In particular, if $q\in(\overline{q},2^{**})$ and $|\alpha|<\overline\alpha$ then global ground states of (\ref{eq:pb}) are not radially symmetric.
\end{Theorem}

The previous result is obtained by noticing that $S_{0,2^{**}}(\lambda)<S_{0,2^{**}}^{\mathrm{rad}}(\lambda)$ and using some continuity of the mappings $(\alpha,q)\mapsto S_{\alpha,q}(\lambda)$ and $(\alpha,q)\mapsto S_{\alpha,q}^{\mathrm{rad}}(\lambda)$. We have no sharp information on the region of values $(\alpha,q)$ for which global ground states of (\ref{eq:pb}) are non radial. On the other hand, for fixed $q$ and $\alpha$, again symmetry breaking is displayed for $\lambda$ large, as stated in the next result.

\begin{Theorem}
\label{T:non-radial-ground-state}
Let $n\ge 5$ and assume (\ref{eq:alpha-beta}). Let $2^{*}=\frac{2n}{n-2}$ be the critical exponent for the first order Sobolev embedding. 
\begin{itemize}
\item[(i)] If $q\in(2^{*},2^{**})$ then for $\lambda$ large enough (depending on $q$) any global ground state of (\ref{eq:pb}) is not radially symmetric.
\item[(ii)] If $q\in(2,2^{*}]$ and 
\begin{equation}
\label{eq:FelliSchneider}
\frac{1}{n-1}\left(\frac{n-4+\alpha}{2}\right)^{2}>\frac{1}{q-2}-\frac{1}{q+2}
\end{equation}
then for $\lambda$ large enough (depending on $q$) any global ground state of (\ref{eq:pb}) is not radially symmetric. 
\end{itemize}
\end{Theorem}

When $q\in(2,2^{*})$ we can better describe the limit profile of ground states of (\ref{eq:pb}) as $\lambda\to\infty$. To this aim, we need to introduce the lower order problem
\begin{equation}
\label{eq:pb1}
\left\{\begin{array}{l}
-\mathrm{div}(|x|^{\alpha-2}\nabla u)=|x|^{-\beta}|u|^{q-2}u\quad\text{in }\R^{n}\\
\int_{\R^{n}}|x|^{\alpha-2}|\nabla u|^{2}~\!dx<\infty
\end{array}\right.
\end{equation}
The natural variational space for problem (\ref{eq:pb1}) is $D^{1,2}(\R^{n};|x|^{\alpha-2})$ defined as the completion of $C^{\infty}_{c}(\R^{n}\setminus\{0\})$ with respect to the norm 
$$
\|u\|_{1,\alpha-2}^{2}=\int_{\R^{n}}|x|^{\alpha-2}|\nabla u|^{2}~\!dx~\!.
$$
Ground states of problem (\ref{eq:pb1}) are defined as weak solutions of (\ref{eq:pb1}) minimizing
$$
\widetilde{S}_{\alpha,q}=\inf_{\scriptstyle u\in D^{1,2}(\R^{n};|x|^{\alpha-2})\atop\scriptstyle u\ne 0}\frac{\int_{\R^{n}}|x|^{\alpha-2}|\nabla u|^{2}~\!dx}{\left(\int_{\R^{n}}|x|^{-\beta}|u|^{q}~\!dx\right)^{2/q}}.
$$
We finally can show:
\begin{Theorem}
\label{T:limiting-profile}
Let $n\ge 5$, $q\in(2,2^{*})$ and assume (\ref{eq:alpha-beta}). If $\lambda_{k}\to\infty$ and $u_{k}\in D^{2,2}(\R^{n};|x|^{\alpha})$ is a ground state of (\ref{eq:pb}) with $\lambda=\lambda_{k}$, then there exists a sequence $(\rho_{k})\subset(0,\infty)$ such that, for a subsequence, $\lambda_{k}^{-\frac{1}{q-2}}\rho_{k}*u_{k}$ converges strongly in $D^{1,2}(\R^{n};|x|^{\alpha-2})$ to a ground state of (\ref{eq:pb1}). The same holds for radial ground states.
\end{Theorem}

We point out that condition (\ref{eq:FelliSchneider}) has been found in \cite{FS} for having symmetry breaking of ground states of the lower order problem (\ref{eq:pb1}). 

We finally observe that most of the results contained in this work are presented in \cite{Ctesi}, generalize in a nontrivial way and complete with new contributions some previous results discussed in \cite{C14} limited to the case $\alpha=0$.

\section{Preliminaries}

Here we introduce the space $D^{2,2}(\R^{n};|x|^{\alpha})$ for $n\ge 5$ and $\alpha\in(4-n,n)$ and we discuss its main embedding properties. The starting point can be given by the weighted Rellich inequality stating that
\begin{equation}
\label{eq:Rellich}
\delta_{\alpha}\int_{\R^{n}}|x|^{\alpha-4}|u|^{2}~\!dx\le\int_{\R^{n}}|x|^{\alpha}|\Delta u|^{2}~\!dx\quad\forall u\in C^{\infty}_{c}(\R^{n}\setminus\{0\})
\end{equation}
with optimal constant
\begin{equation}
\label{eq:deltaalpha}
\delta_{\alpha}=\left[\left(\frac{n-2}{2}\right)^{2}-\left(\frac{\alpha-2}{2}\right)^{2}\right]^2.
\end{equation}
We refer to the paper \cite{CM-CalVar12} and to its bibliography for a deeper discussion on (\ref{eq:Rellich}) and some generalizations. Inequality (\ref{eq:Rellich}) allows us to define the space $D^{2,2}(\R^{n};|x|^{\alpha})$ as the completion of $C^{\infty}_{c}(\R^{n}\setminus\{0\})$ with respect to the Hilbertian norm given by (\ref{eq:norm}). 
Let $2^{**}=\frac{2n}{n-4}$ be the critical exponent for the second order Sobolev embedding. It is known (see, e.g., \cite{CM-Milan11}) that if $q\in[2,2^{**}]$ and $\beta$ is given as in (\ref{eq:alpha-beta}) then $D^{2,2}(\R^{n};|x|^{\alpha})$ is continuously embedded into $L^{q}(\R^{n};|x|^{\beta})$, that is the space of mappings in $L^{q}$ with respect to the measure $|x|^{-\beta}dx$.
\medskip

For future convenience let us introduce also the space $D^{1,2}(\R^{n};|x|^{\widetilde\alpha})$ which can be defined as the completion of $C^{\infty}_{c}(\R^{n}\setminus\{0\})$ with respect to the Hilbertian norm
$$
\|u\|^2_{1,\widetilde\alpha}=\int_{\R^{n}}|x|^{\widetilde\alpha}|\nabla u|^{2}~\!dx~\!.  
$$
This definition of $D^{1,2}(\R^{n};|x|^{\widetilde\alpha})$ is well posed when $\widetilde\alpha>2-n$, thanks to the weighted Hardy inequality
\begin{equation}
\label{eq:Hardy}
h_{\widetilde\alpha}\int_{\R^{n}}|x|^{\widetilde\alpha-2}|u|^{2}~\!dx\le\int_{\R^{n}}|x|^{\widetilde\alpha}|\nabla u|^{2}~\!dx\quad\forall u\in C^{\infty}_{c}(\R^{n}\setminus\{0\})
\end{equation}
which holds with optimal constant
$$
h_{\widetilde\alpha}=\left(\frac{n-2+\widetilde\alpha}{2}\right)^{2}.
$$
In fact here we consider the case $\widetilde\alpha=\alpha-2$. Let $2^{*}=\frac{2n}{n-2}$ be the critical exponent for the first order Sobolev embedding. As a direct consequence of the Caffarelli-Kohn-Nirenberg inequalities \cite{CKN}, if $q\in[2,2^{*}]$ and $\beta$ is given as in (\ref{eq:alpha-beta}) then $D^{1,2}(\R^{n};|x|^{\alpha-2})$ is continuously embedded into $L^{q}(\R^{n};|x|^{\beta})$. 
\medskip

As noticed in the Introduction, for $n\ge 5$ and $\alpha\in(4-n,n)$, the space $D^{2,2}(\R^{n};|x|^{\alpha})$ turns out to be continuously embedded into $D^{1,2}(\R^{n};|x|^{\alpha-2})$. We denote by $\gamma_{\alpha}$ the embedding constant, defined as in (\ref{eq:delta-nabla}). The value of $\gamma_{\alpha}$ is explicitly known only for $\alpha$ in a sub-interval of the admissible one $(4-n,n)$. More precisely, there exists $\alpha^{*}\in(4-n,0)$ such that $\gamma_{\alpha}=(n-\alpha)^{2}/4$ for $\alpha\in(\alpha^{*},n)$, whereas this expression is not valid for all $\alpha\in(4-n,n)$ (see \cite{CM-Bangalore}).  
\medskip

Hence, for $n\ge 5$, $q\in(2,2^{**}]$, $\alpha$ and $\beta$ as in (\ref{eq:alpha-beta}), and $\lambda>-\gamma_{\alpha}$, the minimization problem (\ref{eq:Sq}) is meaningful since the corresponding infimum value $S_{\alpha,q}(\lambda)$ is positive. Let us recall the following well known property linking the minimization problem (\ref{eq:Sq}) with (\ref{eq:pb}).

\begin{Lemma}
\label{L:inf=equation}
If $S_{\alpha,q}(\lambda)>0$ and $u\in D^{2,2}(\R^{n};|x|^{\alpha})$ is a minimum point for (\ref{eq:Sq}), then $U=\tau^{\frac{1}{q-2}}u$ is a solution of (\ref{eq:pb}), being $\tau=S_{\alpha,q}(\lambda)\left(\int_{\R^{n}}|x|^{-\beta}|u|^{q}~\!dx\right)^{\frac{2-q}{q}}$.
\end{Lemma}

Since we are interested also in radial ground states of (\ref{eq:pb}), we introduce also spaces of radial functions. More precisely, with obvious notation, we set
\begin{gather*}
D^{2,2}_{rad}(\R^{n};|x|^{\alpha}):=\{u\in D^{2,2}(\R^{n};|x|^{\alpha})~|~u=u(|x|)\},\\
D^{1,2}_{rad}(\R^{n};|x|^{\widetilde\alpha}):=\{u\in D^{1,2}(\R^{n};|x|^{\widetilde\alpha})~|~u=u(|x|)\},
\end{gather*}
with the same restrictions on $n$, $\alpha$ and $\widetilde\alpha$ as before. For $n\ge 5$ and $\alpha\in(4-n,n)$ the value of the best constant embedding of $D^{2,2}_{rad}(\R^{n};|x|^{\alpha})$ into $D^{1,2}_{rad}(\R^{n};|x|^{\alpha-2})$ is sharp:
\begin{equation}
\label{eq:nabla-delta-rad}
\inf_{\scriptstyle u\in D^{2,2}_{rad}(\R^{n};|x|^{\alpha})\atop\scriptstyle u\ne 0}\frac{\int_{\R^{n}}|x|^{\alpha}|\Delta u|^{2}~\!dx}{\int_{\R^{n}}|x|^{\alpha-2}|\nabla u|^{2}~\!dx}=\frac{(n-\alpha)^{2}}{4}
\end{equation}
(see \cite{CM-Bangalore} and \cite{GM11}).
 
The Emden-Fowler transform, defined by
\begin{equation}
\label{eq:EmdenFowler}
u(x)=|x|^{\frac{4-n-\alpha}{2}}w(-\log|x|)
\end{equation}
provides a nice isomorphism between the space $D^{2,2}_{rad}(\R^{n};|x|^{\alpha})$ and the standard Sobolev space $H^{2}(\R)$, and between $D^{1,2}_{rad}(\R^{n};|x|^{\alpha-2})$ and $H^{1}(\R)$. Indeed:

\begin{Lemma}
\label{L:EmdenFowler}
Let $n\ge 5$ and $\alpha\in(4-n,n)$. For any radial mapping $u\colon\R^{n}\setminus\{0\}\to\R$ let $w\colon\R\to\R$ be defined by (\ref{eq:EmdenFowler}), and viceversa. 
\begin{itemize}
\item[(i)]
$u\in D^{2,2}_{rad}(\R^{n};|x|^{\alpha})$ if and only if $w\in H^{2}(\R)$. Moreover
$$
\int_{\R^{n}}|x|^{\alpha}|\Delta u|^{2}~\!dx=\omega_{n}\int_{\R}\left(|w''|^{2}+2\widetilde{\delta}_{\alpha}|w'|^{2}+\delta_{\alpha}|w|^{2}\right)~\!dt
$$
where $\omega_{n}=|\mathbb{S}^{n-1}|$ and $\widetilde{\delta}_{\alpha}=\left(\frac{n-2}{2}\right)^{2}+\left(\frac{\alpha-2}{2}\right)^{2}$.
\item[(ii)]
$u\in D^{1,2}_{rad}(\R^{n};|x|^{\alpha-2})$ if and only if $w\in H^{1}(\R)$. Moreover
$$
\int_{\R^{n}}|x|^{\alpha-2}|\nabla u|^{2}~\!dx=\omega_{n}\int_{\R}\left(|w'|^{2}+\widetilde{h}_{\alpha}|w|^{2}\right)~\!dt
$$
where $\widetilde{h}_{\alpha}={h}_{\alpha-2}=\left(\frac{n-4+\alpha}{2}\right)^{2}$.
\item[(iii)]
For $q\ge 2$ and $\beta$ as in (\ref{eq:alpha-beta}), $u\in L^{q}(\R^{n};|x|^{-\beta})$ if and only if $w\in L^{q}(\R)$. In this case
$$
\int_{\R^{n}}|x|^{-\beta}|u|^{q}~\!dx=\omega_{n}\int_{\R}|w|^{q}~\!dt~\!.
$$
\end{itemize}
\end{Lemma}
For a proof we refer to \cite{CM-Milan11}.

Thanks to Lemma \ref{L:EmdenFowler} the spaces $D^{2,2}_{rad}(\R^{n};|x|^{\alpha})$ and $D^{1,2}_{rad}(\R^{n};|x|^{\alpha-2})$ are embedded into $L^{q}(\R^{n};|x|^{-\beta})$ for \emph{every} $q\ge 2$. In particular, taking account also of (\ref{eq:nabla-delta-rad}), the minimization problem (\ref{eq:Sqrad}) is meaningful for all $\lambda>-\frac{(n-\alpha)^{2}}{4}$ and $q>2$. Moreover, an analogue of Lemma \ref{L:inf=equation} holds true, namely:

\begin{Lemma}
\label{L:radial-inf=equation}
If $S_{\alpha,q}^{rad}(\lambda)>0$ and $u\in D^{2,2}_{rad}(\R^{n};|x|^{\alpha})$ is a minimum point for (\ref{eq:Sqrad}), then $U=\tau^{\frac{1}{q-2}}u$ is a solution of (\ref{eq:pb}), being $\tau=S_{\alpha,q}^{rad}(\lambda)\left(\int_{\R^{n}}|x|^{-\beta}|u|^{q}~\!dx\right)^{\frac{2-q}{q}}$.
\end{Lemma}

\section{Global ground states}
In this section we prove Theorem \ref{T:ground-state}. 

\begin{Lemma}
\label{L:Ekeland}
There exists a sequence $(u_{k})$ in $D^{2,2}(\R^{n};|x|^{\alpha})$ satisfying
\begin{gather}
\label{eq:lambda-norm}
\int_{\R^{n}}\left(|x|^{\alpha}|\Delta u_{k}|^2+\lambda|x|^{\alpha-2}|\nabla u_{k}|^{2}\right)~\!dx=S_{\alpha,q}(\lambda)^{{q}/{q-2}}+o(1)\\
\label{eq:q-norm}
\int_{\R^{n}}|x|^{-\beta}|u_{k}|^{q}~\!dx=S_{\alpha,q}(\lambda)^{{q}/{q-2}}\\ 
\label{eq:B2}
\int_{B_{2}}|x|^{-\beta}|u_{k}|^{q}~\!dx=\frac{1}{2}S_{\alpha,q}(\lambda)^{{q}/{q-2}}\\
\label{eq:almost-solution}
\Delta(|x|^{\alpha}\Delta u_{k})-\lambda~\!\mathrm{div}(|x|^{\alpha-2}\nabla u_{k})-|x|^{-\beta}|u_{k}|^{q-2}u_{k}\to 0\text{ in $(D^{2,2}(\R^{n};|x|^{\alpha}))'$.}
\end{gather}
\end{Lemma}

\proof
Set 
$$
F(u)=\int_{\R^{n}}|x|^{\alpha}|\Delta u|^{2}~\!dx+\lambda\int_{\R^{n}}|x|^{\alpha-2}|\nabla u|^{2}~\!dx\quad\text{and}\quad G(u)=\int_{\R^{n}}|x|^{-\beta}|u|^{q}~\!dx
$$
so that
$$
S_{\alpha,q}(\lambda)=\inf\{F(u)~|~G(u)=1\}.
$$
Since the constraint $G(u)=1$ defines a smooth manifold, by the Ekeland variational principle, one can find a sequence $(v_{k})\subset D^{2,2}(\R^{n};|x|^{\alpha})$ such that
$$
F'(v_{k})-\mu_{k}G'(v_{k})\to 0\text{~~in~~}(D^{2,2}(\R^{n};|x|^{\alpha}))'\text{~~where~~}\mu_{k}=\frac{\langle\nabla F(v_{k}),\nabla G(v_{k})\rangle}{\|\nabla G(v_{k})\|_{2,\alpha}},
$$
$\nabla F(v_{k})$, $\nabla G(v_{k})$ denote the Riesz representative in $D^{2,2}(\R^{n};|x|^{\alpha})$ of the functionals $F'(v_{k})$, $G'(v_{k})$, respectively, and $\langle~\!,\rangle$ stands for the inner product in $D^{2,2}(\R^{n};|x|^{\alpha})$ corresponding to the norm (\ref{eq:norm}). One can easily check that the sequence $(v_{k})$ is bounded in $D^{2,2}(\R^{n};|x|^{\alpha})$, $\mu_{k}\to\frac{2}{q}S_{\alpha,q}(\lambda)$ and $\sup_{k}\|G'(v_{k})\|<\infty$. Hence
$$
F'(v_{k})-\mu G'(v_{k})\to 0\text{~~in~~}(D^{2,2}(\R^{n};|x|^{\alpha}))'\text{~~where~~}\mu=\frac{2}{q}S_{\alpha,q}(\lambda)~\!.
$$
Now the sequence
$$
\tilde{u}_{k}=S_{\alpha,q}(\lambda)^{1/(q-2)}v_{k}
$$
turns out to satisfy (\ref{eq:lambda-norm}), (\ref{eq:q-norm}) and (\ref{eq:almost-solution}). Finally, for every $k$ one can find $\rho_{k}>0$ such that 
$$
u_{k}(x)=\rho_{k}^{\frac{n-4+\alpha}{2}}\tilde{u}_{k}(\rho_{k} x)
$$
satisfies (\ref{eq:B2}). The sequence $(u_{k})$ always verifies (\ref{eq:lambda-norm}), (\ref{eq:q-norm}) and (\ref{eq:almost-solution}) because the functionals $F$ and $G$ are invariant with respect to (\ref{eq:dilation}).
\QED 

A key tool in our argument is the following compactness lemma. This result is an adaptation of a tool already used in previous works, like \cite{BM} or \cite{CM-Milan11}.

\begin{Lemma}
\label{L:compact}
Let $(u_{k})$ be a sequence in $D^{2,2}(\R^{n};|x|^{\alpha})$ satisfying (\ref{eq:almost-solution}). If
\begin{gather}
\label{eq:weak-conv}
u_{k}\to 0\text{ weakly in }D^{2,2}(\R^{n};|x|^{\alpha})\\
\label{eq:small}
\limsup\int_{B_{R}}|x|^{-\beta}|u_{k}|^{q}~\!dx<S_{\alpha,q}(\lambda)^{q/(q-2)}\text{ for some $R>0$,}
\end{gather}
then $|x|^{-\beta}|u_{k}|^{q}\to 0$ strongly in $L^{1}_{\mathrm{loc}}(B_{R})$.
\end{Lemma}

\proof
Fix $R'\in(0,R)$ and take a cut-off function $\varphi\in C^{\infty}_{c}(B_{R})$ such that $\varphi=1$ on $B_{R'}$. We point out that the sequence $(\varphi^{2} u_{k})$ is bounded in $D^{2,2}(\R^{n};|x|^{\alpha})$. Using $\varphi^{2} u_{k}$ as a test function in (\ref{eq:almost-solution}) we obtain
\begin{equation}
\label{eq:almost-sol-test}
\begin{split}
\int_{\R^{n}}\varphi^{2} u_{k}\Delta(|x|^{\alpha}\Delta u_{k})~\!dx-\lambda&\int_{\R^{n}}\varphi^{2}u_{k}~\!\mathrm{div}(|x|^{\alpha-2}\nabla u_{k})~\!dx\\
&=\int_{\R^{n}}|x|^{-\beta}\varphi^{2}|u_{k}|^{q}~\!dx+o(1).   
\end{split}
\end{equation}
By (\ref{eq:weak-conv}) $u_{k}\to 0$ weakly in $H^{2}_{\mathrm{loc}}(\R^{n}\setminus\{0\})$ and then, by compactness, $u_{k}\to 0$ strongly in $H^{1}_{\mathrm{loc}}(\R^{n}\setminus\{0\})$. Hence we have that
\begin{gather*}
\int_{\R^{n}}|x|^{\alpha}|\Delta(\varphi u_{k})|^{2}~\!dx=\int_{\R^{n}}|x|^{\alpha}\varphi^{2}|\Delta u_{k}|^{2}~\!dx+o(1)\\
\int_{\R^{n}}|x|^{\alpha}(\Delta u_{k})\Delta(\varphi^{2} u_{k})~\!dx=\int_{\R^{n}}|x|^{\alpha}\varphi^{2}|\Delta u_{k}|^{2}~\!dx+o(1)\\
\int_{\R^{n}}|x|^{\alpha-2}\nabla u_{k}\cdot\nabla(\varphi^{2}u_{k})~\!dx=
\int_{\R^{n}}|x|^{\alpha-2}|\nabla(\varphi u_{k})|^{2}~\!dx+o(1).
\end{gather*}
Then, after integration by parts,
\begin{gather*}
\label{eq:Delta2}
\int_{\R^{n}}\varphi^{2} u_{k}\Delta(|x|^{\alpha}\Delta u_{k})~\!dx=\int_{\R^{n}}|x|^{\alpha}|\Delta(\varphi u_{k})|^{2}~\!dx+o(1)\\
\label{eq:nabla}
\int_{\R^{n}}\varphi^{2}u_{k}~\!\mathrm{div}(|x|^{\alpha-2}\nabla u_{k})~\!dx =-\int_{\R^{n}}|x|^{\alpha-2}|\nabla(\varphi u_{k})|^{2}~\!dx+o(1).
\end{gather*}
Consequently (\ref{eq:almost-sol-test}) reduces to
\begin{equation}
\label{eq:almost-sol-test2}
\int_{\R^{n}}|x|^{\alpha}|\Delta(\varphi u_{k})|^{2}~\!dx+\lambda\int_{\R^{n}}|x|^{\alpha-2}|\nabla(\varphi u_{k})|^{2}~\!dx=\int_{\R^{n}}|x|^{-\beta}\varphi^{2}|u_{k}|^{q}~\!dx+o(1).
\end{equation}
By (\ref{eq:small}) there exists $\varepsilon_{0}>0$ such that
\begin{equation}
\label{eq:epsilon}
\int_{B_{R}}|x|^{-\beta}|u_{k}|^{q}~\!dx\le\varepsilon_{0}<S_{\alpha,q}(\lambda)^{q/(q-2)}\quad\forall k\text{ large.}
\end{equation}
Therefore, using the H\"{o}lder inequality and (\ref{eq:epsilon}), we estimate
\begin{equation}
\label{eq:almost-sol-test3}
\int_{\R^{n}}|x|^{-\beta}\varphi^{2}|u_{k}|^{q}~\!dx\le \varepsilon_{0}^{(q-2)/q}\left(\int_{\R^{n}}|x|^{-\beta}|\varphi u_{k}|^{q}~\!dx\right)^{2/q}.
\end{equation}
On the other side, by definition of $S_{\alpha,q}(\lambda)$,
\begin{equation}
\label{eq:almost-sol-test4}
\int_{\R^{n}}|x|^{\alpha}|\Delta(\varphi u_{k})|^{2}~\!dx+\lambda\int_{\R^{n}}|x|^{\alpha-2}|\nabla(\varphi u_{k})|^{2}~\!dx\ge S_{\alpha,q}(\lambda)\left(\int_{\R^{n}}|x|^{-\beta}|\varphi u_{k}|^{q}~\!dx\right)^{\frac{2}{q}}.
\end{equation}
Therefore from (\ref{eq:almost-sol-test2})--(\ref{eq:almost-sol-test4}) it follows that
$$
S_{\alpha,q}(\lambda)\left(\int_{\R^{n}}|x|^{-\beta}|\varphi u_{k}|^{q}~\!dx\right)^{2/q}\le\varepsilon_{0}^{(q-2)/q}\left(\int_{\R^{n}}|x|^{-\beta}|\varphi u_{k}|^{q}~\!dx\right)^{2/q}+o(1).
$$
As $\varepsilon_{0}<S_{\alpha,q}(\lambda)^{q/(q-2)}$ we infer that 
$$
\int_{\R^{n}}|x|^{-\beta}|\varphi u_{k}|^{q}~\!dx\to 0
$$
and then, since $\varphi=1$ on $B_{R'}$ and $R'$ is arbitrary in $(0,R)$, $|x|^{-\beta}|u_{k}|^{q}\to 0$ strongly in $L^{1}_{\mathrm{loc}}(B_{R})$.
\QED
\medskip

\noindent
\textbf{Proof of Theorem \ref{T:ground-state}, part (i).}
Let $(u_{k})$ be a sequence in $D^{2,2}(\R^{n};|x|^{\alpha})$ satisfying (\ref{eq:lambda-norm})--(\ref{eq:almost-solution}), as given by Lemma \ref{L:Ekeland}. Since $\lambda>-\gamma_{\alpha}$, by (\ref{eq:lambda-norm}), the sequence $(u_{k})$ is bounded in $D^{2,2}(\R^{n};|x|^{\alpha})$ and then it admits a subsequence, still denoted $(u_{k})$, weakly converging to some $u\in D^{2,2}(\R^{n};|x|^{\alpha})$. If $u\ne 0$, then $u$ is a minimizer for $S_{\alpha,q}(\lambda)$ and $u_{k}\to u$ strongly in $D^{2,2}(\R^{n};|x|^{\alpha})$. The proof of this fact is definitely standard: one can adapt to our situation a well known argument (see, e.g., \cite{Str}, Chapt.~1, Sect.~4). Hence we have to exclude that $u=0$. We argue by contradiction, assuming that $u=0$. In this case, by Lemma \ref{L:compact},
\begin{equation}
\label{eq:B1-zero}
\int_{B_{1}}|x|^{-\beta}|u_{k}|^{q}\!~dx\to 0.
\end{equation}
Therefore, by (\ref{eq:B2}),
\begin{equation}
\label{eq:annulus}
\int_{B_{2}\setminus B_{1}}|x|^{-\beta}|u_{k}|^{q}\!~dx\to\frac{1}{2}S_{q}(\lambda)^{{q}/{q-2}}.
\end{equation}
Since $q\in(2,2^{**})$ and $u_{k}\to 0$ weakly in $H^{2}_{\mathrm{loc}}(\R^{n}\setminus\{0\})$, the Rellich compactness Theorem implies that $u_{k}\to 0$ strongly in $L^{q}(B_{2}\setminus B_{1})$, in contradiction with (\ref{eq:annulus}). Therefore $u$ cannot be zero and the proof is complete.
\QED
\medskip

Now we focus on the case of critical exponent $q=2^{**}$. To this purpose let us denote by $S^{**}$ the best constant for the second order standard Sobolev embedding, defined by
$$
S^{**}:=\inf_{\scriptstyle u\in C^{\infty}_{c}(\R^{n})\atop\scriptstyle u\ne 0}\frac{\int_{\R^{n}}|\Delta u|^{2}~\!dx}{\left(\int_{\R^{n}}|u|^{2^{**}}~\!dx\right)^{2/2^{**}}}.
$$
As $S^{**}>0$ one can introduce the space $D^{2,2}(\R^{n})$ as the completion of $C^{\infty}_{c}(\R^{n})$ with respect to the norm $\|\Delta u\|_{L^{2}}$. One has that $C^{\infty}_{c}(\R^{n}\setminus\{0\})$ is dense in $D^{2,2}(\R^{n})$ and then $D^{2,2}(\R^{n})$ coincides with the space $D^{2,2}(\R^{n};|x|^{\alpha})$ with $\alpha=0$ and $S^{**}=S_{0,2^{**}}(0)$. Let us recall the following result.

\begin{Lemma}\label{L:EFJ}(\cite{EFJ})
The function $U(x)={\left(1+|x|^{2}\right)^{-\frac{n-4}{2}}}$ is a minimizer  for $S^{**}$ in $D^{2,2}(\R^{n})$.
\end{Lemma}

A condition for existence of a ground state for problem (\ref{eq:pb}) in case of critical exponent is stated by the following result.
\begin{Lemma}
\label{L:<}
If $S_{\alpha,2^{**}}(\lambda)<S^{**}$ then $S_{\alpha,2^{**}}(\lambda)$ is attained in $D^{2,2}(\R^{n};|x|^{\alpha})$.
\end{Lemma}

\proof
As in the proof of Theorem \ref{T:ground-state}, part (i), there exists a sequence $(u_{k})$ in $D^{2,2}(\R^{n};|x|^{\alpha})$ satisfying (\ref{eq:lambda-norm})--(\ref{eq:almost-solution}) and there exists $u\in D^{2,2}(\R^{n};|x|^{\alpha})$ such that $u_{k}\to u$ weakly in $D^{2,2}(\R^{n};|x|^{\alpha})$ and strongly in $H^{1}_{\mathrm{loc}}(\R^{n}\setminus\{0\})$. If $u\ne 0$ then, with a stardard argument, $u$ turns out to be a minimizer. Assume by contradiction that $u=0$. Then (\ref{eq:B1-zero}) and (\ref{eq:annulus}) hold. Let us fix a cut-off function $\varphi\in C^{\infty}_{c}(\R^{n}\setminus\{0\})$ such that $0\le\varphi\le 1$ and $\varphi(x)=1$ for $1\le|x|\le 2$. Arguing as in the first part of the proof of Lemma \ref{L:compact} we obtain (\ref{eq:almost-sol-test2}). By (\ref{eq:lambda-norm}) and (\ref{eq:q-norm}), we also have that
\begin{equation}
\label{eq:almost-sol-crit11}
\begin{split}
\int_{\R^{n}}|x|^{\alpha}|\Delta(\varphi u_{k})|^{2}~\!dx&+\lambda\int_{\R^{n}}|x|^{\alpha-2}|\nabla(\varphi u_{k})|^{2}~\!dx\\
&\le S_{\alpha,2^{**}}(\lambda)\left(\int_{\R^{n}}|x|^{-\beta}|\varphi u_{k}|^{2^{**}}~\!dx\right)^{\frac{2}{2^{**}}}+o(1).
\end{split}
\end{equation}
Since
$$
\int_{\R^{n}}|x|^{\alpha-2}|\nabla(\varphi u_{k})|^{2}~\!dx\le C\int_{B_{2}\setminus B_{1}}\left(|\nabla u_{k}|^{2}+u_{k}^{2}\right)~\!dx
$$
for some constant $C>0$, and $u_{k}\to 0$ strongly in $H^{1}_{\mathrm{loc}}(\R^{n}\setminus\{0\})$, (\ref{eq:almost-sol-test2}) reduces to
\begin{equation}
\label{eq:almost-sol-crit2}
\int_{\R^{n}}|x|^{\alpha}|\Delta(\varphi u_{k})|^{2}~\!dx=\int_{\R^{n}}|x|^{-\beta}\varphi^{2}|u_{k}|^{2^{**}}~\!dx+o(1).
\end{equation}
Now we apply the identity
\begin{equation*}
\begin{split}
|x|^{-\alpha}|\Delta(|x|^{\frac{\alpha}{2}}w)|^{2}&=|\Delta w|^{2}+\alpha^{2}|x|^{-4}|x \cdot \nabla w|^{2}+\frac{\alpha^{2}}{4}\left(n-2+\frac{\alpha}{2}\right)^{2}|x|^{-4}w^{2}\\
&\quad+2\alpha|x|^{-2}(x\cdot\nabla w)\Delta w+\alpha^{2}\left(n-2+\frac{\alpha}{2}\right)|x|^{-4}w(x\cdot\nabla w)\\
&\quad+\alpha\left(n-2+\frac{\alpha}{2}\right)|x|^{-2}w\Delta w
\end{split}
\end{equation*}
with $w=\varphi u_{k}$ and, using again the fact that $u_{k}\to 0$ strongly in $H^{1}_{\mathrm{loc}}(\R^{n}\setminus\{0\})$, we infer that
$$
\int_{\R^{n}}|x|^{\alpha}|\Delta(\varphi u_{k})|^{2}~\!dx=\int_{\R^{n}}\left|\Delta(|x|^{\alpha/2}\varphi u_{k})\right|^{2}~\!dx+o(1).
$$
Hence from (\ref{eq:almost-sol-test2}), (\ref{eq:almost-sol-crit11}) and (\ref{eq:almost-sol-crit2}) it follows that
\begin{equation*}
\begin{split}
S_{\alpha,2^{**}}(\lambda)\left(\int_{\R^{n}}|x|^{-\beta}|\varphi u_{k}|^{2^{**}}~\!dx\right)^{\frac{2}{2^{**}}}&\ge\int_{\R^{n}}\left|\Delta(|x|^{\alpha/2}\varphi u_{k})\right|^{2}~\!dx+o(1)\\
&\ge S^{**}\left(\int_{\R^{n}}||x|^{\alpha/2}\varphi u_{k}|^{2^{**}}~\!dx\right)^{\frac{2}{2^{**}}}+o(1).
\end{split}
\end{equation*}
Since 
$$
\int_{\R^{n}}||x|^{\alpha/2}\varphi u_{k}|^{2^{**}}~\!dx=\int_{\R^{n}}|x|^{-\beta}|\varphi u_{k}|^{2^{**}}~\!dx
$$
and, by hypothesis, $S_{\alpha,2^{**}}(\lambda)<S^{**}$, we deduce that
$$
\int_{\R^{n}}|x|^{-\beta}|\varphi u_{k}|^{2^{**}}~\!dx\to 0
$$
in contradiction with (\ref{eq:annulus}), as $\varphi=1$ on $B_{2}\setminus B_{1}$.
\QED

\begin{Lemma}
\label{L:lambda}
If 
\begin{equation}
\label{eq:condition-lambda}
-\gamma_{\alpha}<\lambda<\left(\alpha-\frac{\alpha^{2}}{4}\right)\frac{(\alpha^{2}-4\alpha)(n-3)-2(n-2)^{2}(n-4)}{(n-4+\alpha)^{2}(n-3)+(n-4)^{2}}
\end{equation}
then $S_{\alpha,2^{**}}(\lambda)<S^{**}$.
\end{Lemma}

\proof
Set $u(x)=|x|^{-\frac{\alpha}{2}}U(x)$ where $U$ is as in Lemma \ref{L:EFJ}. One can check that
\begin{gather*}
\int_{\R^{n}}|x|^{\alpha}|\Delta u|^{2}~\!dx=\int_{\R^{n}}|\Delta U|^{2}~\!dx+A_{\alpha}\int_{\R^{n}}|x|^{-4}U^{2}~\!dx
\\
\int_{\R^{n}}|x|^{\alpha-2}|\nabla u|^{2}~\!dx=B_{\alpha}\int_{\R^{n}}|x|^{-4}U^{2}~\!dx
\\
\int_{\R^{n}}|x|^{-\beta}|u|^{2^{**}}~\!dx=\int_{\R^{n}}U^{2^{**}}~\!dx
\end{gather*}
where
\begin{gather*}
A_{\alpha}=\left(\frac{\alpha^{2}}{4}-\alpha\right)\left(\frac{\alpha^{2}}{4}-\alpha-\frac{(n-2)^{2}(n-4)}{2(n-3)}\right)\\
B_{\alpha}=\left(\frac{n-4+\alpha}{2}\right)^{2}+\frac{(n-4)^{2}}{4(n-3)}.
\end{gather*}
For detailed computation see \cite{CM-Milan11} or \cite{Ctesi}.
Then
$$
S_{\alpha,q}(\lambda)\le\frac{\int_{\R^{n}}|\Delta U|^{2}~\!dx}{\left(\int_{\R^{n}}U^{2^{**}}~\!dx\right)^{2/2^{**}}}+\left(A_{\alpha}+\lambda~\!B_{\alpha}\right)\frac{\int_{\R^{n}}U^{2}~\!dx}{\left(\int_{\R^{n}}U^{2^{**}}~\!dx\right)^{2/2^{**}}}~\!.
$$
Hence the strict inequality $S_{\alpha,q}(\lambda)<S^{**}$ holds true when
\begin{equation}
\label{eq:lambda-n-alpha}
\lambda<\lambda_{\alpha}:=-\frac{A_{\alpha}}{B_{\alpha}}=\left(\alpha-\frac{\alpha^{2}}{4}\right)\frac{(\alpha^{2}-4\alpha)(n-3)-2(n-2)^{2}(n-4)}{(n-4+\alpha)^{2}(n-3)+(n-4)^{2}}.
\end{equation}
We point out that $\lambda_{\alpha}=0$ when $\alpha=0,4$, $\lambda_{\alpha}>0>-\gamma_{\alpha}$ when $\alpha\in(4-n,0)\cup(4,n)$, whereas $\lambda_{\alpha}<0$ when $\alpha\in(0,4)$. In fact, with some calculation one can check that in this last case $\lambda_{\alpha}>-\frac{(n-\alpha)^{2}}{4}=-\gamma_{\alpha}$. Hence (\ref{eq:condition-lambda}) is completely proved.
\QED

Clearly the proof of Theorem \ref{T:ground-state}, part (ii) is a direct consequence of Lemmas \ref{L:<} and \ref{L:lambda}. 

When $\alpha=0$, condition (\ref{eq:condition-lambda}) is optimal for the validity of the strict inequality $S_{\alpha,2^{**}}(\lambda)<S^{**}$. Indeed one has:

\begin{Proposition}
\label{P}
If $\lambda\ge 0$ then $S_{0,2^{**}}(\lambda)=S^{**}$. Moreover, for $\lambda>0$ the infimum $S_{0,2^{**}}(\lambda)$ is not achieved in $D^{2,2}(\R^{n})$. 
\end{Proposition}
For a proof, see \cite{C14} or \cite{Ctesi}.

\section{Radial ground states}
Here we prove Theorem \ref{T:radial-ground-state}. In view of Lemma \ref{L:EmdenFowler} we have that
\begin{equation}
\label{eq:radial-minimization}
S_{\alpha,q}^{\mathrm{rad}}(\lambda)=\omega_{n}^{\frac{q-2}{q}}\inf_{\scriptstyle w\in H^{2}(\R)\atop\scriptstyle w\ne 0}\frac{\int_{\R}\left(|w''|^{2}+2a_{\lambda}|w'|^{2}+b_{\lambda}|w|^{2}\right)~\!dt}{\left(\int_{\R}|w|^{q}~\!dt\right)^{{2}/{q}}}
\end{equation}
where
\begin{equation}
\label{eq:alambda-blambda}
a_{\lambda}=\frac{(n-2)^{2}}{4}+\frac{(\alpha-2)^{2}}{4}+\frac{\lambda}{2}\quad\text{and}\quad b_{\lambda}=\left(\frac{(n-\alpha)^{2}}{4}+\lambda\right)\left(\frac{n-4+\alpha}{2}\right)^{2}.
\end{equation}
We point out that, thanks to the assumption $\lambda>-(n-\alpha)^{2}/{4}$, the values $a_{\lambda}$ and $b_{\lambda}$ are positive. Now we use the following key result, proved in \cite{BM}:

\begin{Theorem}
\label{T:BM}
For every $a,b>0$ and $q>2$ the minimization problem
$$
\inf_{\scriptstyle w\in H^{2}(\R)\atop\scriptstyle w\ne 0}\frac{\int_{\R}\left(|w''|^{2}+2a|w'|^{2}+b|w|^{2}\right)~\!dt}{\left(\int_{\R}|w|^{q}~\!dt\right)^{{2}/{q}}}
$$
admits a minimum point. In addition, if $a^{2}\ge b$ then the minimum point is positive and unique, up to the natural invariances of the problem (i.e., translation, inversion, multiplication by a non zero constant).
\end{Theorem}

In the case in consideration 
$$
a_{\lambda}^{2}-b_{\lambda}=\left(\frac{(n-2)(\alpha-2)}{2}-\frac{\lambda}{2}\right)^{2}.
$$
Hence by Theorem \ref{T:BM} there exists a positive function $w\in H^{2}(\R)$ which is a minimizer for the problem defined by the right hand side of (\ref{eq:radial-minimization}). Such a minimizer is unique up to translation, inversion, and multiplication by a non zero constant. Then, using Lemma \ref{L:EmdenFowler}, we infer that the mapping $u$ defined by (\ref{eq:EmdenFowler}) belongs to $D^{2,2}_{rad}(\R^{n};|x|^{\alpha})$, is a positive minimizer for $S_{\alpha,q}^{rad}(\lambda)$ and is the unique minimizer up to the weighted dilation (\ref{eq:dilation}) and to a multiplicative constant. Then one applies Lemma \ref{L:radial-inf=equation} to get a radial ground state for problem (\ref{eq:pb}).
\QED

\section{Symmetry breaking and limiting profiles}
This section contains the proof of the symmetry breaking results stated in Theorems \ref{T:SB1} and \ref{T:non-radial-ground-state}, and the description of the limit profile of ground states for $q\in(2,2^{*})$, when $\lambda\to\infty$ (Theorem \ref{T:limiting-profile}). 

Let us start with the discussion of Theorem \ref{T:SB1}, whose proof lies on the following semicontinuity inequalities.

\begin{Lemma}
\label{L:semicontinuity-inequalities}
For $\lambda\ge 0$ one has
\begin{equation}
\label{eq:sc-ineq}
\limsup_{(\alpha,q)\to(0,2^{**}_{-})}S_{\alpha,q}(\lambda)\le S_{0,2^{**}}(\lambda)\quad\text{and}\quad S^{\mathrm{rad}}_{0,2^{**}}(\lambda)\le\liminf_{(\alpha,q)\to(0,2^{**}_{-})}S^{\mathrm{rad}}_{\alpha,q}(\lambda)~\!.
\end{equation}
\end{Lemma}

\proof
Fix $\lambda\ge 0$. For every $u\in C^{\infty}_{c}(\R^{n}\setminus\{0\})$, $u\ne 0$, set
$$
Q_{\alpha,q}(u)=\frac{\int_{\R^{n}}|x|^{\alpha}|\Delta u|^{2}~\!dx+\lambda\int_{\R^{n}}|x|^{\alpha-2}|\nabla u|^{2}~\!dx}{\left(\int_{\R^{n}}|x|^{-\beta}|u|^{q}~\!dx\right)^{{2}/{q}}}.
$$
Since $Q_{\alpha,q}(u)\to Q_{0,2^{**}}(u)$ as $(\alpha,q)\to(0,2^{**}_{-})$ and $C^{\infty}_{c}(\R^{n}\setminus\{0\})$ is dense in $D^{2,2}(\R^{n};|x|^{\alpha})$, the first inequality in (\ref{eq:sc-ineq}) immediately follows. In order to check the second inequality, we proceed in this way. For every $q\in(2,2^{**})$ let $\theta_{q}=\frac{q-2}{2^{**}-2}$. By the H\"{o}lder inequality, one has that
$$
\int_{\R^{n}}|x|^{-\beta}|u|^{q}~\!dx\le\left(\int_{\R^{n}}|x|^{\alpha-4}|u|^{2}~\!dx\right)^{1-\theta_{q}}\left(\int_{\R^{n}}|x|^{\frac{n\alpha}{n-4}}|u|^{2^{**}}~\!dx\right)^{\theta_{q}}
$$
and 
consequently
\begin{equation}
\label{eq:sn1}
S_{\alpha,q}^{\mathrm{rad}}(\lambda)^{\frac{q}{\theta_{q}2^{**}}}\ge S_{\alpha,2}^{\mathrm{rad}}(\lambda)^{\frac{(1-\theta_{q})2}{\theta_{q}2^{**}}}S_{\alpha,2^{**}}^{\mathrm{rad}}(\lambda)\ge \delta_{\alpha}^{\frac{(1-\theta_{q})2}{\theta_{q}2^{**}}}S_{\alpha,2^{**}}^{\mathrm{rad}}(\lambda)
\end{equation}
(see (\ref{eq:Rellich})--(\ref{eq:deltaalpha})).
Now we estimate $S_{\alpha,2^{**}}^{\mathrm{rad}}(\lambda)$ in terms of $S_{0,2^{**}}^{\mathrm{rad}}(\lambda)$. 
For every $\alpha$ let $\tau_{\alpha}:=1+\frac{\alpha}{n-4}$ and for every radial $u\in C^{\infty}_{c}(\R^{n}\setminus\{0\})$, set $\widetilde{u}(x)=u(|x|^{1/\tau_{\alpha}})$. One can check that
\begin{gather*}
\int_{\R^{n}}|\widetilde{u}|^{2^{**}}~\!dx=\tau_{\alpha}\int_{\R^{n}}\!|x|^{\frac{n\alpha}{n-4}}|u|^{2^{**}}~\!dx~\!,~~\int_{\R^{n}}\!|x|^{-2}|\nabla\widetilde{u}|^{2}~\!dx=\tau_{\alpha}^{-1}\!\int_{\R^{n}}\!|x|^{\alpha-2}|\nabla u|^{2}~\!dx,\\
\int_{\R^{n}}\!|\Delta\widetilde{u}|^{2}~\!dx=\tau_{\alpha}^{-3}\!\int_{\R^{n}}\!|x|^{\alpha}|\Delta u+R_{\alpha}u|^{2}~\!dx~~\text{where}~ R_{\alpha}u=(\tau_{\alpha}-1)(n-2)\nabla u\cdot\frac{x}{|x|^{2}}.
\end{gather*}
Setting $\varepsilon_{\alpha}=|\tau_{\alpha}-1|(n-2)\gamma_{\alpha}^{-1/2}$ and using (\ref{eq:delta-nabla}), as $\alpha\in(4-n,n)$ we can estimate
$$
\int_{\R^{n}}|x|^{\alpha}|R_{\alpha}u|^{2}~\!dx\le\varepsilon_{\alpha}^{2}\int_{\R^{n}}|x|^{\alpha}|\Delta u|^{2}~\!dx
$$
and then, by the Cauchy-Schwarz inequality,
$$
\int_{\R^{n}}|x|^{\alpha}|\Delta u+R_{\alpha}u|^{2}~\!dx\le (1+\varepsilon_{\alpha})^{2}\int_{\R^{n}}|x|^{\alpha}|\Delta u|^{2}~\!dx~\!.
$$
Therefore
\begin{equation*}
\begin{split}
S_{0,2^{**}}^{\mathrm{rad}}(\lambda)&\le\frac{\int_{\R^{n}}(|\Delta\widetilde{u}|^{2}+\lambda|x|^{-2}|\nabla\widetilde{u}|^{2})~\!dx}{\left(\int_{\R^{n}}|\widetilde{u}|^{2^{**}}~\!dx\right)^{\frac{2}{2^{**}}}}\\
&\le\frac{(1+\varepsilon_{\alpha})^{2}}{\tau_{\alpha}^{3+2/2^{**}}}~\!
\frac{\int_{\R^{n}}\left(|x|^{\alpha}|\Delta u|^{2}+\lambda|x|^{\alpha-2}|\nabla u|^{2}\right)~\!dx}{\left(\int_{\R^{n}}|x|^{\frac{n\alpha}{n-4}}|u|^{2^{**}}~\!dx\right)^{\frac{2}{2^{**}}}}\\
&\quad+\frac{\lambda}{\tau_{\alpha}^{1+2/2^{**}}}\left(1-\frac{(1+\varepsilon_{\alpha})^{2}}{\tau_{\alpha}^{2}}\right)\frac{\int_{\R^{n}}|x|^{\alpha-2}|\nabla u|^{2}~\!dx}{\left(\int_{\R^{n}}|x|^{\frac{n\alpha}{n-4}}|u|^{2^{**}}~\!dx\right)^{\frac{2}{2^{**}}}}\\
&\le K_{\alpha}\frac{\int_{\R^{n}}\left(|x|^{\alpha}|\Delta u|^{2}+\lambda|x|^{\alpha-2}|\nabla u|^{2}\right)~\!dx}{\left(\int_{\R^{n}}|x|^{\frac{n\alpha}{n-4}}|u|^{2^{**}}~\!dx\right)^{\frac{2}{2^{**}}}}
\end{split}
\end{equation*}
where
$$
K_{\alpha}={\left(\frac{(1+\varepsilon_{\alpha})^{2}}{\tau_{\alpha}^{3+2/2^{**}}}+\frac{\lambda\gamma_{\alpha}^{-1}}{\tau_{\alpha}^{1+2/2^{**}}}\left|1-\frac{(1+\varepsilon_{\alpha})^{2}}{\tau_{\alpha}^{2}}\right|\right)}.
$$
Hence
\begin{equation}
\label{eq:sn2}
S_{\alpha,2^{**}}^{\mathrm{rad}}(\lambda)\ge K_{\alpha}^{-1}S_{0,2^{**}}^{\mathrm{rad}}(\lambda)~\!.
\end{equation}
For $|\alpha|$ small, one has that $C^{-1}\le\gamma_{\alpha}\le C$ and $C^{-1}\le\delta_{\alpha}\le C$ for some constant $C>0$. Consequently, if $(\alpha,q)\to(0,2^{**}_{-})$, then $\tau_{\alpha}\to 1$, $\theta_{q}\to 1$, $\varepsilon_{\alpha}\to 0$, and $K_{\alpha}\to 1$. Thus the second inequality in (\ref{eq:sc-ineq}) follows from (\ref{eq:sn1})--(\ref{eq:sn2}). 
\QED
\medskip

\noindent
\textbf{Proof of Theorem \ref{T:SB1}.} Fix $\lambda>0$. By Theorem \ref{T:radial-ground-state} $S_{0,2^{**}}^{\mathrm{rad}}(\lambda)$ is achieved in $D^{2,2}_{rad}(\R^{n})$. Instead, by Proposition \ref{P}, $S_{0,2^{**}}(\lambda)$ is not attained. Hence $S_{0,2^{**}}^{\mathrm{rad}}(\lambda)>S_{0,2^{**}}(\lambda)$. Therefore the conclusion follows by applying Lemma \ref{L:semicontinuity-inequalities}.
\QED

Now let us address to Theorem \ref{T:non-radial-ground-state}. It is convenient to normalize the minimization problems defined by (\ref{eq:Sq}) and (\ref{eq:Sqrad}) as follows. For every $\varepsilon\ge 0$ set
\begin{gather*}
\widetilde{S}_{\alpha,q}(\varepsilon):=\inf_{\scriptstyle u\in C^{\infty}_{c}(\R^{n}\setminus\{0\})\atop\scriptstyle u\ne 0}\frac{\varepsilon\int_{\R^{n}}|x|^{\alpha}|\Delta u|^{2}~\!dx+\int_{\R^{n}}|x|^{\alpha-2}|\nabla u|^{2}~\!dx}{\left(\int_{\R^{n}}|x|^{-\beta}|u|^{q}~\!dx\right)^{{2}/{q}}}\\
\widetilde{S}_{\alpha,q}^{rad}(\varepsilon):=\inf_{\scriptstyle u\in C^{\infty}_{c}(\R^{n}\setminus\{0\})\atop\scriptstyle u=u(|x|),~u\ne 0}\frac{\varepsilon\int_{\R^{n}}|x|^{\alpha}|\Delta u|^{2}~\!dx+\int_{\R^{n}}|x|^{\alpha-2}|\nabla u|^{2}~\!dx}{\left(\int_{\R^{n}}|x|^{-\beta}|u|^{q}~\!dx\right)^{{2}/{q}}}
\end{gather*}
We remark that if $\lambda>0$ then
\begin{equation}
\label{eq:lambda-eps}
\widetilde{S}_{\alpha,q}(\lambda^{-1})=\lambda^{-1}{S}_{\alpha,q}(\lambda)\quad\text{and}\quad\widetilde{S}_{\alpha,q}^{rad}(\lambda^{-1})=\lambda^{-1}{S}_{\alpha,q}^{rad}(\lambda).
\end{equation}

\begin{Lemma}\label{L:Seps}
Let $n\ge 5$ and $\alpha\in(4-n,n)$. 
\begin{itemize}
\item[(i)]
For every $q\in(2,2^{**}]$ one has that $\widetilde{S}_{\alpha,q}(\varepsilon)\to\widetilde{S}_{\alpha,q}(0)$  as $\varepsilon\to 0$.
\item[(ii)] 
For every $q>2$ one has that $\widetilde{S}_{\alpha,q}^{rad}(\varepsilon)\to\widetilde{S}_{\alpha,q}^{rad}(0)$ as $\varepsilon\to 0$.
\end{itemize}
\end{Lemma}

\proof (i) Fix $\alpha\in(4-n,n)$ and $q\in(2,2^{**}]$ and set
$$
\widetilde{Q}_{\varepsilon}(u):=\frac{\varepsilon\int_{\R^{n}}|x|^{\alpha}|\Delta u|^{2}~\!dx+\int_{\R^{n}}|x|^{\alpha-2}|\nabla u|^{2}~\!dx}{\left(\int_{\R^{n}}|x|^{-\beta}|u|^{q}~\!dx\right)^{{2}/{q}}}\quad\forall u\in C^{\infty}_{c}(\R^{n}\setminus\{0\}),~u\ne 0.
$$
Since $\widetilde{S}_{\alpha,q}(\varepsilon)\le\widetilde{Q}_{\varepsilon}(u)$, when $\varepsilon\to 0$ one has that
$$
\limsup \widetilde{S}_{\alpha,q}(\varepsilon)\le\widetilde{Q}_{0}(u)\quad\forall u\in C^{\infty}_{c}(\R^{n}\setminus\{0\}),~u\ne 0
$$
namely $\limsup\widetilde{S}_{\alpha,q}(\varepsilon)\le\widetilde{S}_{\alpha,q}(0)$. On the other hand $\widetilde{Q}_{\varepsilon}(u)\ge \widetilde{Q}_{0}(u)$ and then $\widetilde{S}_{\alpha,q}(\varepsilon)\ge\widetilde{S}_{\alpha,q}(0)$. Hence the conclusion immediately follows. Clearly (ii) is proved in the same way.
\QED
\medskip

\noindent
\textbf{Proof of Theorem \ref{T:non-radial-ground-state}.}
If $q\in(2^{*},2^{**}]$ then $\widetilde{S}_{\alpha,q}^{rad}(0)>0$ whereas $\widetilde{S}_{\alpha,q}(0)=0$. Indeed, by Lemma \ref{L:EmdenFowler}, 
$$
\widetilde{S}_{\alpha,q}^{rad}(0)=\omega_{n}^{1-\frac{2}{q}}\inf_{\scriptstyle w\in C^{\infty}_{c}(\R)\atop\scriptstyle w\ne 0}\frac{\int_{\R}(|w'|^{2}+\widetilde{h}_{\alpha}|w|^{2})~\!dt}{\left(\int_{\R}|w|^{q}~\!dt\right)^{2/q}}
$$
which is positive because $H^{1}(\R)$ is embedded into $L^{q}$. Instead, taking $u\in C^{\infty}_{c}(\R^{n})$, $u\ne 0$, with support in the unit ball, fixing $x_{0}\in\R^{n}$ with $|x_{0}|=1$ and setting
$$
u_{\delta}(x)=\delta^{-\frac{n-2}{2}}u\left(\frac{x-x_{0}}{\delta}\right),
$$
one can check that $\widetilde{Q}_{0}(u_{\delta})\to 0$ as $\delta\to 0$, because $q>2^{*}$. Hence, by Lemma \ref{L:Seps}, $\widetilde{S}_{\alpha,q}(\varepsilon)<\widetilde{S}_{\alpha,q}^{rad}(\varepsilon)$ for $\varepsilon>0$ small, and then ${S}_{\alpha,q}(\lambda)<{S}_{\alpha,q}^{rad}(\lambda)$ for $\lambda$ large, by (\ref{eq:lambda-eps}). If $q\in(2,2^{*}]$ and (\ref{eq:FelliSchneider}) holds, then $\widetilde{S}_{\alpha,q}(0)<\widetilde{S}_{\alpha,q}^{rad}(0)$, as proved in \cite{FS} and one concludes as before that ${S}_{\alpha,q}(\lambda)<{S}_{\alpha,q}^{rad}(\lambda)$ for $\lambda$ large.
\QED

In the following we study the behavior of ground states of problems (\ref{eq:pb}) for fixed $q\in(2,2^{*})$, in the limit $\lambda\to\infty$.

\begin{Lemma}
\label{L:limit-profile}
Let $n\ge 5$, $q\in(2,2^{*})$ and assume (\ref{eq:alpha-beta}). Let $\varepsilon_{k}\to 0^{+}$ and for every $k$ let $v_{k}\in D^{2,2}(\R^{n};|x|^{\alpha})$ be a minimizer for $\widetilde{S}_{\alpha,q}(\varepsilon_{k})$ with 
\begin{equation}
\label{eq:assum-k}
\int_{\R^{n}}|x|^{-\beta}|v_{k}|^{q}~\!dx=1\quad\text{and}\quad\int_{B_{2}}|x|^{-\beta}|v_{k}|^{q}~\!dx=\frac{1}{2}~\!.
\end{equation}
If $v_{k}\to v$ weakly in $D^{1,2}(\R^{n};|x|^{\alpha-2})$ then $v_{k}\to v$ strongly in $D^{1,2}(\R^{n};|x|^{\alpha-2})$ and $v$ is a minimizer for $\widetilde{S}_{\alpha,q}(0)$. 
\end{Lemma}

\proof
Let us write, briefly, $\widetilde{S}_{k}=\widetilde{S}_{\alpha,q}(\varepsilon_{k})$ and $\widetilde{S}_{0}=\widetilde{S}_{\alpha,q}(0)$. Since $v_{k}$ is a minimizer for $\widetilde{S}_{k}$ and (\ref{eq:assum-k}) holds, we have that
\begin{equation}
\label{eq:tilde-bounded}
\widetilde{S}_{0}\le\int_{\R^{n}}|x|^{\alpha-2}|\nabla v_{k}|^{2}~\!dx\le\widetilde{S}_{k}
\end{equation}
and then, 
\begin{equation}
\label{eq:zero-k}
\varepsilon_{k}\int_{\R^{n}}|x|^{\alpha}|\Delta v_{k}|^{2}~\!dx\to 0
\end{equation}
for $\widetilde{S}_{k}\to\widetilde{S}_{0}$ by Lemma \ref{L:Seps}. Now we want to exclude that $v=0$. To do this, we argue by contradiction, assuming that $v\to 0$ weakly in $D^{1,2}(\R^{n};|x|^{\alpha-2})$. Since $v_{k}$ is a minimizer for $\widetilde{S}_{k}$, it is so for $S_{\alpha,q}(\varepsilon_{k}^{-1})$ and, by Lemma \ref{L:inf=equation}, 
\begin{equation}
\label{eq:pb-epsk}
\varepsilon_{k}\Delta(|x|^{\alpha}\Delta v_{k})-\mathrm{div}(|x|^{\alpha-2}\nabla v_{k})=\widetilde{S}_{k}|x|^{-\beta}|v_{k}|^{q-2}v_{k}\quad\text{on }\R^{n}.
\end{equation}
Taking a cut-off function $\varphi\in C^{\infty}_{c}(\R^{n})$ with $\mathrm{supp}(\varphi)\subset B_{2}$ and $\varphi\equiv 1$ in $B_{1}$, we can use $\varphi^{2}v_{k}$ as a test function in (\ref{eq:pb-epsk}) getting that
\begin{equation}
\label{eq:start-k}
\varepsilon_{k}\int_{\R^{n}}|x|^{\alpha}\Delta v_{k}\Delta(\varphi^{2}v_{k})~\!dx+\int_{\R^{n}}|x|^{\alpha-2}\nabla v_{k}\cdot\nabla(\varphi^{2}v_{k})~\!dx=\widetilde{S}_{k}\int_{\R^{n}}|x|^{-\beta}\varphi^{2}|v_{k}|^{q}~\!dx~\!.
\end{equation}
We estimate each term of (\ref{eq:start-k}) as follows. Firstly
\begin{equation*}
\begin{split}
&\left|\int_{\R^{n}}|x|^{\alpha}\Delta v_{k}\Delta(\varphi^{2}v_{k})~\!dx\right|\le\int_{\R^{n}}|x|^{\alpha}|\Delta v_{k}|\big|\Delta(\varphi^{2})v_{k}+2\nabla(\varphi^{2})\cdot\nabla v_{k}\big|~\!dx\\
&\hspace{7.8cm}+\int_{\R^{n}}|x|^{\alpha}|\varphi\Delta v_{k}|^{2}~\!dx\\
&\le\left(\int_{\R^{n}}|x|^{\alpha}|\Delta v_{k}|^{2}~\!dx\right)^{1/2}\left(C\int_{B_{2}\setminus B_{1}}(|v_{k}|^{2}+|\nabla v_{k}|^{2})~\!dx\right)^{1/2}+\int_{\R^{n}}|x|^{\alpha}|\Delta v_{k}|^{2}~\!dx
\end{split}
\end{equation*}
because $\varphi$ is constant in $B_{1}$ and out of $B_{2}$. Then
\begin{equation}
\label{eq:Delta-k}
\varepsilon_{k}\int_{\R^{n}}|x|^{\alpha}\Delta v_{k}\Delta(\varphi^{2}v_{k})~\!dx\to 0
\end{equation}
thanks to (\ref{eq:zero-k}) and because the sequence $(v_{k})$ is bounded in $D^{1,2}(\R^{n};|x|^{\alpha-2})$ and then also in $H^{1}(B_{2}\setminus B_{1})$. 
Secondly, as $v_{k}\to 0$ weakly in $D^{1,2}(\R^{n};|x|^{\alpha-2})$ we also obtain that 
\begin{equation}
\label{eq:grad-k}
\int_{\R^{n}}|x|^{\alpha-2}\nabla v_{k}\cdot\nabla(\varphi^{2}v_{k})~\!dx=\int_{\R^{n}}|x|^{\alpha-2}|\nabla(\varphi v_{k})|^{2}~\!dx+o(1)~\!.
\end{equation}
In addition, using H\"older inequality and (\ref{eq:assum-k})
\begin{equation}
\label{eq:nonlinear-k}
\int_{\R^{n}}|x|^{-\beta}|v_{k}|^{q}\varphi^{2}~\!dx\le 2^{-\frac{q-2}{q}}\left(\int_{\R^{n}}|x|^{-\beta}|\varphi v_{k}|^{q}~\!dx\right)^{2/q}.
\end{equation}
By definition of $\widetilde{S}_{0}$ we also have that
\begin{equation}
\label{eq:CKN-k}
\int_{\R^{n}}|x|^{\alpha-2}|\nabla(\varphi v_{k})|^{2}~\!dx\ge\widetilde{S}_{0}\left(\int_{\R^{n}}|x|^{-\beta}|\varphi v_{k}|^{q}~\!dx\right)^{2/q}.
\end{equation}
Plugging (\ref{eq:Delta-k})--(\ref{eq:CKN-k}) into (\ref{eq:start-k}) and taking into account that $\widetilde{S}_{k}\to\widetilde{S}_{0}$, we obtain
$$
\widetilde{S}_{0}\left(\int_{\R^{n}}|x|^{-\beta}|\varphi v_{k}|^{q}~\!dx\right)^{2/q}\le 2^{-\frac{q-2}{q}}\widetilde{S}_{0}\left(\int_{\R^{n}}|x|^{-\beta}|\varphi v_{k}|^{q}~\!dx\right)^{2/q}+o(1)~\!.
$$
Being $q>2$, we infer that $\int_{B_{1}}|x|^{-\beta}|v_{k}|^{q}~\!dx\to 0$ and, by (\ref{eq:assum-k}), 
\begin{equation}
\label{eq:k-Lq-annulus}
\int_{B_{2}\setminus B_{1}}|x|^{-\beta}|v_{k}|^{q}~\!dx\to\frac{1}{2}~\!.
\end{equation}
On the other hand, if $v_{k}\to 0$ weakly in $D^{1,2}(\R^{n};|x|^{\alpha-2})$, in particular $v_{k}\to 0$ weakly in $H^{1}_{loc}(\R^{n}\setminus\{0\})$ and, by the Rellich compactness theorem, $v_{k}\to  0$ strongly in $L^{q}(B_{2}\setminus B_{1})$ because $q<2^{*}$. This is in contradiction with (\ref{eq:k-Lq-annulus}). Therefore $v\ne 0$. Since, by (\ref{eq:zero-k}), $(v_{k})$ is a minimizing sequence for $\widetilde{S}_{0}$ with a nonzero weak limit, it is standard to conclude that $v$ is a minimizer for $\widetilde{S}_{0}$ and $v_{k}\to v$ strongly in $D^{1,2}(\R^{n};|x|^{\alpha-2})$.
\QED
\medskip

\noindent
\textbf{Proof of Theorem \ref{T:limiting-profile}.}
Let $u_{k}\in D^{2,2}(\R^{n};|x|^{\alpha})$ be a ground state of (\ref{eq:pb}) with $\lambda=\lambda_{k}$ and set $S_{k}=S_{\alpha,q}(\lambda_{k})$ and $\widetilde{u}_{k}=S_{k}^{-\frac{1}{q-2}}u_{k}$. Then $\widetilde{u}_{k}$ turns out to be a minimizer for $\widetilde{S}_{\alpha,q}(\lambda_{k}^{-1})$ with
$$
\int_{\R^{n}}|x|^{-\beta}|\widetilde{u}_{k}|^{q}~\!dx=1~\!.
$$
Moreover there exists $\rho_{k}>0$ such that 
$$
\int_{B_{2}}|x|^{-\beta}|\rho_{k}*\widetilde{u}_{k}|^{q}~\!dx=\frac{1}{2}~\!.
$$
Hence $v_{k}=\rho_{k}*\widetilde{u}_{k}$ is again a minimizer for $S_{\alpha,q}(\lambda_{k}^{-1})$ and satisfies (\ref{eq:assum-k}). Since (\ref{eq:tilde-bounded}) holds, there exists $v\in D^{1,2}(\R^{n};|x|^{\alpha})$ such that, for a subsequence, $v_{k}\to v$ weakly in $D^{1,2}(\R^{n};|x|^{\alpha})$. Since $\lambda_{k}\to\infty$, by Lemma \ref{L:limit-profile}, $v_{k}\to v$ strongly in $D^{1,2}(\R^{n};|x|^{\alpha})$ and $v$ is a minimizer for $\widetilde{S}_{\alpha,q}(0)$. Then $u=\widetilde{S}_{\alpha,q}(0)^{\frac{1}{q-2}}v$ turns out to be a ground state for problem (\ref{eq:pb1}) and $\lim_{k}\lambda_{k}^{-\frac{1}{q-2}}\rho_{k}*u_{k}=\widetilde{S}_{\alpha,q}(0)^{\frac{1}{q-2}}\lim_{k} S_{k}^{-\frac{1}{q-2}}\rho_{k}*u_{k}=u$ strongly in $D^{1,2}(\R^{n};|x|^{\alpha})$.
\QED
\medskip

\noindent
\textbf{Acknowledgments.} The first author is partially supported by the PRIN2012 grant ``Variational and perturbative aspects of nonlinear differential problems''. He is also member of the Gruppo Nazionale per l'Analisi Matematica, la Probabilit\`a e le loro Applicazioni (GNAMPA) of the Istituto Nazionale di Alta Matematica (INdAM).

\label{References}

\end{document}